\documentclass[12pt]{amsart}
\usepackage{amsmath,amsthm,amssymb}

\newcommand{\al}{\alpha}

\newcommand{\fr}{\mathcal{F}}
\newcommand{\ity}{\infty}
\newcommand{\C}{\mathbb{C}}

\newcommand{\N}{\mathbb{N}}

\numberwithin{equation}{section}
\newtheorem{theorem}{Theorem}[section]
\newtheorem{lemma}[theorem]{Lemma}

\newtheorem*{ta}{Theorem A}
\newtheorem*{tb}{Theorem B}
\newtheorem*{tc}{Theorem C}
\newtheorem*{td}{Theorem D}
\newtheorem*{te}{Theorem E}
\newtheorem*{SFT}{Second Fundamental Theorem}
\theoremstyle{remark}

\newtheorem*{exam}{Example}

\makeatletter
\@namedef{subjclassname@2010}{%
  \textup{2010} Mathematics Subject Classification}
\makeatother

\begin{document}

\title[NORMALITY CRITERIA FOR A FAMILY OF MEROMORPHIC FUNCTIONS ]{NORMALITY CRITERIA FOR A FAMILY OF MEROMORPHIC FUNCTIONS WITH MULTIPLE ZEROS}


\author[G. Datt]{Gopal Datt}
\address{School of Mathematics, Harish-Chandra Research Institute,
Chhatnag Road, Jhunsi, Allahabad -- 211019, India} \email{ggopal.datt@gmail.com }\email{gopaldatt@hri.res.in}

\author[Y. T. Li]{Yuntong Li}
\address{Department of Basic Courses, Shaanxi Railway Institute, Weinan 714000, Shaanxi Province, PR China}\email{liyuntong2005@sohu.com}

\author[P. Rani]{Poonam Rani}

\address{Department of Mathematics,University of Delhi,
Delhi -- 110 007, India }
\email{pnmrani753@gmail.com}

\begin{abstract}

In this article, we prove some normality criteria for a family of meromorphic functions having zeros with some multiplicity.   Our main result  involves sharing of a  holomorphic function  by certain  differential polynomials. Our results generalize some of the results of Fang and Zalcman ~\cite{FZ} and Chen et al \cite{CTZY} to a great extent.
\end{abstract}



 \maketitle

\section{Introduction and main results}
One important aspect  of the theory of complex analytic functions is to find normality criteria for families of meromorphic functions. The notion of normal families was introduced by Paul Montel in 1907. Let us begin by recalling the definition. {\it A family of meromorphic (holomorphic) functions defined on a domain $D\subset \C$ is said  to be \emph{normal} in the domain, if every sequence in the family  has a subsequence which converges spherically  uniformly on compact subsets of  $D$ to a meromorphic (holomorphic) function or to $\ity$} \cite{Ahl, Hay, Schiff, Yang}.\\

In ~\cite{MS}  Mues and Steinmetz proved a uniqueness theorem which says that: {\it{if a non-constant meromorphic  function  $f$ in the plane,  shares three distinct complex numbers $a_1, a_2,  a_3$ with its first order derivative $f'$, then $f \equiv f'$.}}
Wilhelm Schwick ~\cite{Sch} was the first who gave a connection between normality and shared values and  proved a theorem related to above result of ~\cite{MS}  which says that:  {\it{the family $\fr$ of meromorphic functions on a domain ${D}$ is normal, if $f$ and $f'$ share $a_1$, $a_2$, $a_3$ for every $f\in \mathcal F$, where $a_1$, $a_2$, $a_3$ are distinct complex numbers.}}

Let us recall the meaning of shared value. Let $f$ be a meromorphic function of a domain $D\subset\C$. For $p\in \C$, Let
\begin{equation*}
  E_f(p)=\{z\in D: f(z)=p\}
\end{equation*}
and let
\begin{equation*}
E_f(\infty)= \text{poles of\ } f \text{\ in\ } D.
\end{equation*}
 For $p\in \C\cup \{\infty\}$, two meromorphic function $f$ and $g$ of $D$ share the value $p$ if $E_f(p)=E_g(p).$\\

In 2008 Fang and Zalcman ~\cite{FZ} proved the following normality criteria:\\

\begin{ta}~\cite{FZ}
Let $\fr$ be a family of meromorphic functions on a domain $D$, let $n\geq 2$ be a positive integer, and $a(\neq 0), b\in \C$. If for each $f\in\fr$, all zeros of $f$ are multiple and $f+a(f')^n\neq b$ on $D$, then $\fr$ is normal on $D$.
\end{ta}
Related to the above result of ~\cite{FZ}, Wang \cite{YMW} proved the following result on normality and sharing value:
\begin{tb}\cite{YMW}
Let $\fr$ be a family of meromorphic functions on the plane domain $D$, let $n\geq3$ be a positive integer. Let $a, b$ be two finite complex numbers such that $a\neq 0$. If  all zeros of $f$ are multiple for each $f\in\fr$, and $f+a(f')^n$ and $g+a(g')^n$ share $b$ in $D$ for every pair of functions $f, g\in \fr$, then $\fr$ is normal in $D$.
\end{tb}
Extending the result of \cite{FZ}, Xu, Wu and Liao ~\cite{XWL} proved the following normality criteria:
\begin{tc}\label{thm B}~\cite{XWL}
Let $\fr$ be a family of meromorphic functions on a plane domain $D$, let  $a(\neq 0), b\in \C$, and $n, k$ be two positive integers such that $n\geq k+1$. If for each $f\in\fr$, $f$ has only zeros of multiplicity at least $k+1$, and $f+a\left(f^{(k)}\right)^n\neq b$ on $D$, then $\fr$ is normal on $D$.
\end{tc}
Related to the result of \cite{XWL}  chen et al ~\cite{CTZY} proved the following normality criteria concerning shared values:
\begin{td}\label{thm A}~\cite{CTZY}
Let $\fr$ be a family of meromorphic functions on the plane domain $D$, let $n, k$ be  positive integers such that $n\geq k+2$, and $a, b$ be two finite complex numbers such that $a\neq 0$. If  all zeros of $f$ have multiplicity at least $k+1$ for each $f\in\fr$, and $f+a\left(f^{(k)}\right)^n$ and $g+a\left(g^{(k)}\right)^n$ share $b$ in $D$ for every pair of functions $f, g\in \fr$, then $\fr$ is normal in $D$.
\end{td}
It is evident from the following questions that the above theorem has not been stated in full generality: \begin{enumerate}
                                    \item[$Q. 1.$]\ Can we  weaken the condition on $n$?
                                    \item[$Q. 2.$]\ Can we replace value $b$ by any holomorphic function?
                                  \end{enumerate}    In this paper, we try to give the  answers of these questions and see that after weakening the condition on $n\geq{k+2}$ to $n>2$ the theorem is valid for the case where multiplicities of zeros of $f\in\fr$ are at least $2k+1$, and $b$ is a non-vanishing holomorphic function.  Now we state our main result.

\begin{theorem}\label{thm 1'}
Let $\alpha(z)\not\equiv 0$ be a holomorphic function  with zeros of multiplicity at most $m$ in $D$. Let $a\in\C$ be a non-zero constant, and  $n, k$ be  positive integers such that $n>k+1$, and $m<k$. Let $\fr$ be a family of meromorphic functions in the domain $D$. Suppose that for each $f\in\fr$,   all zeros of $f$ have multiplicity at least $2k+2$ and all poles (if exist) have multiplicity at least $2k+3$. If for each pair $f, g$ in $\fr$, $f(z)+a\left(f^{(k)}\right)^n(z)$ and $g(z)+a\left(g^{(k)}\right)^n(z)$ share $\alpha(z)$ IM in $D$, then $\fr$ is normal in $D$.
 \end{theorem}

      When $\alpha(z)$ is non-vanishing holomorphic function in $D$, we get the following strengthened result:                            
\begin{theorem}\label{thm 1}
Let $\alpha(z)\neq 0$ be a holomorphic function  with zeros of multiplicity at most $m$ in $D$. Let $a(z)$ be a non-vanishing holomorphic function in $D$ and  $n, k$ be  positive integers such that $n>2$. Let $\fr$ be a family of meromorphic functions in the domain $D$. Suppose that for each $f\in\fr$,   all zeros of $f$ have multiplicity at least $2k+1$. If for each pair $f, g$ in $\fr$, $f(z)+a(z)\left(f^{(k)}\right)^n(z)$ and $g(z)+a(z)\left(g^{(k)}\right)^n(z)$ share $\alpha(z)$ IM in $D$, then $\fr$ is normal in $D$.
 \end{theorem}
We give the following example in support of Theorem \ref{thm 1}.
\begin{exam}
  Let $D=\{z\in\C:0<|z|<1\}$, $n=3$ and $k=1$. Consider the family $\fr=\{jz^3: j\in \N \}$ and $a(z)=1/z^3$, $\al(z)=z^3$. Clearly $\fr$ satisfies all the conditions of $\fr$ and $\fr$ is normal in $D$.
\end{exam}
We also improved Theorem C in the following manner:
\begin{theorem}\label{thm 2}
Let $\mathcal{F}$ be a family of meromorphic functions on a domain ${D}$, let $n, k$ be  positive integers such that $n>2$. Let $b$ be a non-zero finite complex number and $a(z)$ be a non-vanishing holomorphic function. If for each  function $f\in\mathcal{F}$, all zeros of $f$ have multiplicity at least $2k+1$, and $f(z)+a(z)\left(f^{(k)}\right)^n(z)-b$ has at most one zero IM in ${D}$,  then $\mathcal{F}$ is normal in $D$.
\end{theorem}

Also We  prove a theorem on the value distribution of a transcendental meromorphic function.  The following  theorem on value distribution of a zero-free transcendental meromorphic function  is due to  Li \cite{YLI} ( also see \cite{LFZ}).
\begin{te}\label{lemma 5}
 Let $f$ be a transcendental meromorphic function with $f\neq 0$, let $a$ be  non-zero finite complex number, and let $n\geq 2$ and $k$ be two positive integers. Then $f+a\left(f^{(k)} \right)^n$ assumes each value $b\in \C$ infinitely often.
 \end{te}
 In the following theorem we prove above result for  the case where $f\not\equiv 0$.
 \begin{theorem}\label{lemma 6}
 Let $f$ be a transcendental meromorphic function, let $a$ be  non-zero finite complex number, and let $n\geq 3$  be a positive integer. Then $f+a\left(f^{(k)} \right)^n$ assumes each value $b\in \C$ infinitely often.
\end{theorem}

\section{Some Notations and results of Nevanlinna theory}
 Let $\Delta=\{z: |z|<1\}$ be the unit disk and $\Delta(z_0, r):=\{z: |z-z_0|<r\}.$ We use the following standard functions of value distribution theory, namely
\begin{center}
$T(r,f),  m(r,f),  N(r,f)\ \text{and}\ \overline{N}(r,f)$.
\end{center}
We let $S(r,f)$ be any function satisfying
\begin{center}
$S(r,f)=o\big(T(r,f)\big)$,  as $r\rightarrow +\ity,$
\end{center}
 possibly outside of a set with finite measure.\\

  \begin{SFT}Suppose $f(z)$ is meromorphic in the finite plane and non-degenerate into a constant. If $a_{\nu}(\nu= 1, 2,\ldots, q)$ are $q(\geq 3)$ distinct complex numbers (one of them may be infinity), then
 \begin{equation}
 (q-2)T(r, f)\leq \sum_{\nu=1}^{q}\overline{N}\left( r, \frac{1}{f-a_\nu}\right)+S(r,f).
 \end{equation} \end{SFT}

\section{Some Lemmas}
 In order to prove our results we need  the following Lemmas. The well known Zalcman Lemma is a very important tool in the study of normal families. The following is a new version due to  Zalcman ~\cite{Zalc} (also see  \cite{Zalc 1}, p. $814$).
\begin{lemma}\label{lemma 1} Let $\mathcal F$ be a family of meromorphic  functions in the unit disk  $\Delta$, with the property that for every function $f\in \mathcal F,$  the zeros of $f$ are of multiplicity at least $l$ and the poles of $f$ are of multiplicity at least $k$. If $\mathcal F$ is not normal at $z_0$ in $\Delta$, then for  $-l< \alpha <k$, there exist
\begin{enumerate}
\item{ a sequence of complex numbers $z_n \rightarrow z_0$, $|z_n|<r<1$},
\item{ a sequence of functions $f_n\in \mathcal F$, }
\item{ a sequence of positive numbers $\rho_n \rightarrow 0$},
\end{enumerate}
such that $g_n(\zeta)=\rho_n^{\alpha}f_n(z_n+\rho_n\zeta) $ converges to a non-constant meromorphic function $g$ on $\C$ with $g^{\#}(\zeta)\leq g^{\#}(0)=1$. Moreover, $g$ is of order at most two. Here, $g^{\#}(z)=\frac{|g'(z)|}{1+|g(z)|^2}$ is the spherical derivative of $g$.
\end{lemma}
Let $f$ be a non-constant meromorphic function in $\C$. A  differential polynomial $P$ of $f$ is defined by $\displaystyle {P(z):= \sum_{i=1}^{n}\al_{i}(z)\prod_{j=0}^{p}\left(f^{(j)}\left(z\right)\right)^{S_{ij}}},$ where $S_{ij}$'s are non-negative integers and    $\al_i(z)\not\equiv 0$ are small functions of $f$, that is $T(r,\al_i)=o\big(T(r,f)\big)$. The lower degree of the differential polynomial $P$ is defined by $$d(P):= \min_{1\leq i \leq n}\sum_{j=0}^{p}S_{ij} .$$

The following result was proved by Dethloff et al. in ~\cite{dethloff}.
\begin{lemma}\label{lemma1}
Let $a_1, \ldots, a_q$ be distinct non-zero complex numbers.  Let $f$ be a non-constant  meromorphic function and let $P$  be a non-constant  differential polynomial of $f$ with $d(P)\geq 2.$ Then
\begin{equation}\notag
T(r,f)\leq \left(\frac{q\theta(P)+1}{qd(P)-1}\right)\overline{N}\left(r,\frac{1}{f}\right)+\frac{1}{qd(P)-1}\sum_{j=1}^{q}\overline{N}\left(r,\frac{1}{P-a_j}\right)+ S\left(r,f\right),
\end{equation}
for all $r\in [1,+\ity)$ excluding a set of finite Lebesgue measure, where $\displaystyle\theta(P):= \max_{1\leq i \leq n}\sum_{j=0}^{p}jS_{ij}.$\\

 Moreover, in the case of an entire function, we have
 \begin{equation}\notag
T(r,f)\leq \bigg(\frac{q\theta(P)+1}{qd(P)}\bigg)\overline{N}\left(r,\frac{1}{f}\right)+\frac{1}{qd(P)}\sum_{j=1}^{q}\overline{N}\left(r,\frac{1}{P-a_j}\right)+ S(r,f),
\end{equation}
for all $r\in [1,+\ity)$ excluding a set of finite Lebesgue measure.
\end{lemma}
This result was proved  by Hinchliffe in ~\cite{hin} for $q=1$.\\

We now prove some lemmas to establish our results in the next section.
\begin{lemma}\label{lemma 2}
Let $f(z)$ be a transcendental meromorphic function in $\C$. If all zeros of $f(z)$ has multiplicity at least $2k+1$, then for a positive integer $n>2$ , $\left(f^{(k)}\right)^n$ assumes every non-zero finite value $b$ infinitely often.
\end{lemma}
\begin{proof}
Suppose on the contrary that $\left(f^{(k)}\right)^n$ assumes the value $b$ only finitely many times. Then
\begin{equation}\label{eq}
N\left(r, \frac{1}{(\left(f^{(k)}\right)^n-b}\right) = O\left(\log r\right)= S(r, f).
\end{equation}

Without loss of generality we may assume $b=1$. Let $P=\left(f^{(k)}\right)^n$.  It is easy to see that
$$d(P) = n  \ \text{and}\ \theta(P)=nk.$$
 Clearly $d(P)>2.$ So by Lemma \ref{lemma1} we get
\begin{equation}\notag
T(r,f)\leq \bigg(\frac{nk + 1}{n - 1}\bigg)\overline{N}\left(r,\frac{1}{f}\right)+\bigg(\frac{1}{n -  1}\bigg)\overline{N}\left(r,\frac{1}{P-1}\right)+ S(r,f),
\end{equation}
and this gives
\begin{equation}\notag
T(r,f)\leq \bigg(\frac{nk+1}{(n-1)(2k+1)}\bigg)N\left(r,f\right)+ S(r,f),
\end{equation}
and so
\begin{equation}\notag
\bigg(\frac{(n-2)(k+1)}{(n-1)(2k+1)}\bigg)T(r,f)\leq  S(r,f).
\end{equation}
  But this is a contradiction and hence establishes the lemma.
\end{proof}
\begin{lemma}\label{lemma 3}
Let $f(z)$ be a non-constant rational function in $\C$ and $n>2$ be a positive integer. If all zeros of $f(z)$ has multiplicity at least $2k+1$ then $\left(f^{(k)}\right)^n$ has at least two distinct $b$-points, where $b$ is a non-zero complex number.
\end{lemma}
\begin{proof}
On the contrary, assume that $\left(f^{(k)}\right)^n$ has at most one $b$-point. Now there are two cases to consider.\\

\underline{Case 1.} Let $\left(f^{(k)}\right)^n-b=0$ has exactly one zero and it is at $z_0$.\\

First we assume that $f$ is a non-constant polynomial. Set $\left(f^{(k)}(z)\right)^n-b=A(z-z_0)^l$, where  $A$ is a non-zero constant and $l$ is a positive integer such that $l\geq 2(k+1)$. Then $\left(\left(f^{(k)}(z)\right)^n\right)'=Al(z-z_0)^{l-1}.$ This shows that $z_0$ is the only zero of $\left(\left(f^{(k)}(z)\right)^n\right)'$. Since zeros of $\left(f^{(k)}(z)\right)^n$ are multiple, we deduce that $z_0$ is a zero of $\left(f^{(k)}(z)\right)^n$, which is a contradiction, since $b\neq 0$. \\

 Now suppose that $f$ is a non-polynomial rational function with  zeros of multiplicity at least $2k+1$. Clearly, $f^{(k)}(z)$ is non-constat. Let $b_1, b_2,\ldots, b_n$ be $n$ distinct zeros of $w^n=b$. Let $w=f^{(k)}(z)$ then we obtain $\left(f^{(k)}-b_1\right)\left(f^{(k)}-b_2\right)\ldots \left(f^{(k)}-b_n\right)=0.$ Since $z_0$ is a zero of $\left(f^{(k)}\right)^n-b=0$, so for one $j\in\{1, 2,\ldots, n\},$ $f^{(k)}(z_0)=b_j$ and $f^{(k)}(z_0)\neq b_i$ for $i(\neq j)\in \{1, 2,\ldots, n\}$. Thus we have
\begin{equation}\label{eq1}
f^{(k)}(z)=b_j+\frac{A(z-z_0)^l}{Q(z)}\equiv b_i+\frac{B}{Q(z)},
\end{equation}
where $A, B$  are non-zero constants. By \eqref{eq1} we obtain
\begin{equation} \label{eq2}
(b_i-b_j)Q(z)+B=A(z-z_0)^l.
\end{equation}
From \eqref{eq2} we get $l\geq k+1$, since zeros of $Q(z)$ are of multiplicity at least $k+1$. Again from \eqref{eq2} we obtain $Q(z_0)\neq 0.$ After differentiating \eqref{eq1} we have
\begin{equation}\label{eq3}
f^{(k+1)}(z)=\frac{A(z-z_0)^{l-1}Q(z)-A(z-z_0)^lQ'(z)}{Q^2(z)}\equiv \frac{-BQ'(z)}{Q^2(z)},
\end{equation}
which gives, $A(z-z_0)^{l-1}\left(Q(z)-(z-z_0)Q'(z)\right)=-BQ'(z).$ Since $Q(z)$ is a polynomial of degree $l\geq k+1$ whose zeros are other than $z_0$. This shows that $Q(z)-(z-z_0)Q'(z)$ is a non-constant polynomial. Thus we observe that the degree of $Q'(z)$ is at least $l$. This is a contradiction to the fact that the degree of $Q(z)$ is $l$.\\

\underline{Case 2.} Let $\left(f^{(k)}\right)^n\neq b$.  Let $b_1, b_2,\ldots, b_n$ be $n(\geq 3)$ distinct solutions of $w^n=b$.  By Nevanlinna's second fundamental theorem,
\begin{align*}
T\left(r,f^{(k)}\right)\leq \sum_{\nu =1}^{n}\overline{N}\left(r, \frac{1}{f^{(k)}-b_\nu} \right)+ S\left(r, f^{(k)}\right).
\end{align*}
  It follows that $T(r,f^{(k)})=S(r,f^{(k)}),$ which is a contradiction. This completes the proof of Lemma.
\end{proof}
\begin{lemma}\label{lem}\cite{XWL}
Let $n\geq 2$, $k$ be positive integers, let $p$ be a non-zero constant and let $P(z)$ be a polynomial. Then the solution of the differential equation $p\left(W^{(k)}(z)\right)^n+W(z)=P(z)$ must be a polynomial.
\end{lemma}
\begin{lemma}\label{ lemma 1P}
Let $f $ be a transcedental meromorphic function on the complex plane $\C$, let $a(\neq 0)$ be a complex number and let $ n,m,k$ be three positive integers such that $n\geq k+1$ and  $m<k.$
\begin{enumerate}
\item{If $n\geq k+2,$ then}
\begin{equation}\notag
(n-1)T(r,f^{(k)})\leq (k^2+k+1)\overline{N}(r,f)+(k+1)^2 N\left( r,\frac{1}{f+a(f^{(k)})^n-z^m}\right)+S(r,f^{(k)}).
\end{equation}
\item{If $n=k+1$, then}
\begin{equation}\notag
kT(r,f^{(k)})\leq (k^2+1)\overline{N}(r,f)+(k^2+k) N\left( r,\frac{1}{f+a(f^{(k)})^n-z^m}\right)+S(r,f^{(k)}).
\end{equation}

\end{enumerate}
\end{lemma}
\begin{proof}
Let 

\begin{equation}\label{Ann eq q1}
g=f+a(f^{(k)})-z^m
\end{equation}
and
\begin{equation}\label{Ann eq q2}
h=\frac{g^{(k)}}{g}.
\end{equation}
Then $ h\not\equiv 0$. Otherwise ,if $h\equiv 0$, then $g^{(k)}\equiv 0,$ so we conclude that $g$ is a polynomial with degree at most $k-1$. Noting that $n=k+1 \geq 2$, we conclude from \eqref{Ann eq q1} that $f$ must be a polynomial, which is a contradiction. By simple calculation we have
\begin{align*}
g^{(k)}=f^{(k)}+a[(f^{(k})]^{(k)}=f^{(k)}(1+P(f^{(k)})),
\end{align*}
where
\begin{align}
P(f^{(k)})=a\left(f^{(k)}\right)&^{n-k-1}\notag\\
& \times\left(\frac{n!}{(n-k)!}\left(f^{(k+1})\right)^k+\cdots+n\left(f^{(k)}\right)^{k-1}f^{(2k)}\right)\notag
\end{align}
and $P(f^{(k})$ is a homogeneous differential equation in $f^{(k)}$ of degree $n-1$.Then
\begin{equation}\label{Ann eq q3}
gh=f^{(k)}(1+P(f^{(k)}).
\end{equation}It follows from \eqref{Ann eq q1} that $T(r,g)\leq O(T(r,f))$, and so $S(r,g)=S(r,f)$. This and \eqref{Ann eq q2} give 
\begin{equation}\label{Ann eq q4}
m(r,h)=S(r,f)
\end {equation}
using\eqref{Ann eq q2}-\eqref{Ann eq q4} and Nevanlina's first fundamental theorem,we obtain
\begin{align*}
N\left(r,\frac{1}{f^{(k)}}\right)+N\left(r,\frac{1}{P(f^{(k)})+1}\right)&\leq N\left(r,\frac{1}{h}\right)+N\left(r,\frac{1}{g}\right)\notag\\
&\leq N(r,h)+ N\left (r,\frac{1}{g}\right)+S(r,f)\notag\\
&\leq kN(r,f)+(k+1)N\left(r,\frac{1}{g}\right)+S(r,f).\notag\\
\end{align*} 
On the other hand , by Nevanlinna's first fundamental theorem, we get
\begin{align}
m\left(r,\frac{1}{(f^{(k)})^{n-1}}\right)&+m\left(r,\frac{1}{(P(f^{(k)})+1}\right)\notag\\
&\leq m\left(r,\frac{P(f^{(k)})}{(f^{(k)})^{n-1}}\right)+m\left(r,\frac{1}{P(f^{(k)})}\right)+m\left(r,\frac{1}{(P(f^{(k)})+1}\right)\notag\\&\leq m\left(r,\frac{1}{(P(f^{(k)})}\right)+m\left(r,\frac{1}{(P(f^{(k)})+1}\right)+S(r,f)\notag\\
&\leq m\left(r,\frac{1}{P(f^{(k)})}+\frac{1}{(P(f^{(k)})+1}\right)+S(r,f)\label{eq1.1}\\
&\leq m\left(r,\frac{1}{[P(f^{(k)})]'}\right)+m\left(r,\frac{[P(f{(k)})]'}{P(f^{(k)})}+\frac{[P(f{(k)})]'+1}{(P(f^{(k)})+1}\right)+S(r,f)\notag\\
&\leq m\left(r,\frac{1}{[P(f^{(k)})]'}\right)+S(r,f)\notag\\
& =T(r,[P(f^{(k)})]')-N\left(r,\frac{1}{[P(f^{(k)})]'}\right)+S(r,f).\notag
\end{align}
We deduce from \eqref{eq1.1} and Nevanlina's first fundamental theorem that 
\begin{align}
(n-1)T(r,f^{(k)})&\leq N(r,f)+(n-1)N\left(r,\frac{1}{f^{(k)}}\right)\notag\\
&+N\left(r,\frac{1}{P(f^{(k)}+1}\right)-N\left(r,\frac{1}{[P(f^{(k)})]'}\right)+S(r,f).\label{eq1.2}
\end{align}
If $n=k+1$, from\eqref{eq1.1} and \eqref{eq1.2} we have
\begin{equation}\notag
kT(r,f^{(k)})\leq(k^2+1)N(r,f)+(k^2+k)N\left(r,\frac{1}{g})\right)+S(r,f).
\end{equation}
This prove second part of the Lemma. If $n\geq k+2$ and supposing that $z_0$ is a zero of $f^{(k)}$ of multiplicity $l$, we see that $z_0$ is a zero of $[P(f^{(k)})]'$ of multiplicity at least $(n-1)l-k-1.$ This gives 
\begin{align}
(n-1)N\left(r,\frac{1}{f^{(k)}}\right)+N&\left(r,\frac{1}{P(f^{(k)})+1}\right)-N\left(r,\frac{1}{P[(f^{(k)}]'}\right)\notag\\
&\leq(k+1)\overline{N}\left(r,\frac{1}{f^{(k)}}\right)+\overline{N}\left(r,\frac{1}{P(f^{(k)})+1}\right)\label{eq1.3}
\end{align}
Sustituting \eqref{eq1.3} in \eqref{eq1.2}, we have
\begin{equation}\notag
(n-1)T(r,f^{(k)})\leq \overline{N}(r,f)+(k+1)\overline{N}\left(r,\frac{1}{f^{(k)}}\right)+\overline{N}\left(r,\frac{1}{P(f^{(k)})+1}\right).
\end{equation}
which together with leads to
\begin{equation}\notag
(n-1)T(r,f^{(k)})\leq(k^2+k+1)\overline{N}(r,f)+(k+1)^2 N\left(r,\frac{1}{g}\right)+S(r,f).
\end{equation}
This completes the proof of the Lemma.
\end{proof}
\begin{lemma}\label{lemma trans 1p}
Let $f$ be a transcendental meromorphic function on the complex palne $\C$, let $a$ be a non-zero finite complex number and let $n$,$k$and $m$ be three positive integers such that $n\geq k+1$ and $m<k$ then $f+af^{(k)}-z^m$ assumes infinitely many zeros.  
\end{lemma}
\begin{proof}
Suppose that $f+a(f)^{(k)}-z^m$ has finitely many zeros. Since $f$ is transcendental, we have
\begin{equation}\label{Ann eq ar1}
N\left(r,\frac{1}{f+a(f^{(k)})^n-z^m}\right)=O(\log r)=S(r,f).
\end{equation}
If $n=k+1$. It follows from \eqref{Ann eq ar1} and first part of  Lemma \ref{ lemma 1P}     that
\begin{align*}
kT(r,f^{(k)})&\leq(k^2+1)\overline{N}(r,f)+S(r,f)\notag\\
& \leq\frac{k^2+1}{k+1}N(r,f^{(k)})+S(r,f)\notag\\
&\leq\frac{k^2+1}{k+1}T(r,f^{(k)})+S(r,f),\notag
\end{align*}
that is,\\
\begin{equation}\notag
\frac{k-1}{k+1}T(r,f^{(k)})\leq S(r,f^{(k)}).
\end{equation}
This contradicts the fact that $f$ is transcendental.
If $n\geq k+2,$ then using \eqref{Ann eq ar1} and second part of Lemma \ref{ lemma 1P}, we obtain 
\begin{align*}
(k+1)T(r,f^{(k)})&\leq (n-1)T(r,f^{(k)})\notag\\
&\leq(k^2+k+1)\overline{N}(r,f)+S(r,f)\notag\\
&\leq\frac{k^2+k+1}{k+1}N(r,f^{(k)})+S(r,f)\notag\\
&\leq\frac{k^2+k+1}{k+1}T(r,f^{(k)})+S(r,f)\notag
\end{align*}
Then $ T(r,f^{(k)})\leq S(r,f^{(k)})$. But this is impossible since $f$ is transcendental. The theorem is proved.
\end{proof}

\begin{lemma}\label{lemma 2P}
Let $f$ be a non-constant rational function and let $n,m,k$ be three positive integers such that $n>2$ and $m<k$ . Suppose that every zero of $f$ has multiplicity at least $2k+2$ and every pole(if exists) of $f$ has multiplicity at least $2k+3$.  Then $f+\left(f^{(k)}\right)^n-z^m$ has  at least two distinct zeros.
\end{lemma}
\begin{proof} Let us assume that $D(f)(z)-z^m:=f+\left(f^{(k)}\right)^n-z^m$ has atmost one zero. Now we consider the following cases:\\

\underline{Case 1.} $D(f)(z)-z^m$ has exactly one zero $z_0$ with multiplicity $l.$ \\

\underline{Case1.1.} Suppose that $f$ is a non constant polynomial, then we set
\begin{equation}\label{Ann eq b2}
f(z)=A(z-\al_1)^{m_1}\ldots (z-\al_s)^{m_s},
\end{equation}
where $A$ is a non-zero constant, $m_i \geq 2k+2$ are integers.
Now differentiating \eqref{Ann eq b2} $k$-times, we get
\begin{equation}\label{Ann eq b3}
f^{(k)}(z)=(z-\al_1)^{m_1-k}\ldots (z-\al_s)^{m_s-k} h_1(z),
\end{equation}
where $h_1$ is a non zero polynomial with  $\operatorname{deg}(h_1)\leq k(s-1)$
From \eqref{Ann eq b2} and \eqref{Ann eq b3}, we see that
\begin{align}
D(f)(z)-z^m=f(z)+&(f^{(k)})^n(z)=A(z-\al_1)^{m_1}\ldots (z-\al_s)^{m_s}\notag\\
&+(z-\al_1)^{n(m_1-k)}\ldots (z-\al_s)^{n(m_s-k)} h_1^n(z)-z^m\label{eq2.1P}
\end{align}
Now, differentiating \eqref{eq2.1P} $m+1$-times we get \\
\begin{align}
(D(f)&(z))^{(m+1)}\notag\\
&=(z-\al_1)^{m_1-m-1}\ldots (z-\al_s)^{m_s-m-1}[1+(z-\al_1)^{(n-1)m_1-nk}\ldots (z-\al_s)^{(n-1)m_s-nk}g_1(z)]\notag\\
&=(z-\al_1)^{m_1-m-1}\ldots (z-\al_s)^{m_s-m-1}h_2(z)\label{eq2.2P}
\end{align}
where $h_2$ is a non zero constant polynomial.\\

As $D(f)(z)-z^m$ has only one zero then from eq \eqref{eq2.2P},  we get a contradiction.\\

\underline{Case1.2.} Suppose that $f(z)$ is a non-polynomial rational function defined as
\begin{equation}\label{Ann eq a1}
     f(z)=A\frac{(z-\al_1)^{m_1}\ldots (z-\al_s)^{m_s}}{(z-\beta_1)^{n_1}\ldots (z-\beta_t)^{n_t}},
   \end{equation}
where $A$ is a non-zero constant, $m_i\geq 2k+2\ (i=1, 2, \ldots, s)$ and $n_j\geq 2k+3\ (j=1, 2, \ldots, t).$\\
   Let us define
   \begin{equation}\label{Ann eq a2}
     \sum_{i=1}^{s}m_i=M\geq(2k+2)s \ \text{and}\ \sum_{j=1}^{t}n_j=N\geq (2k+3)t.
        \end{equation}
from \eqref{Ann eq a1} it follows that
\begin{equation}\label{Ann eq a3}
     f^{(k)}(z)=A\frac{(z-\al_1)^{m_1-k}\ldots (z-\al_s)^{m_s-k}}{(z-\beta_1)^{n_1+k}\ldots (z-\beta_t)^{n_t+k}}g_1(z),
   \end{equation} 
where $g_1$ is a polynomial with $\operatorname{deg}(g_1)\leq k(s+t-1)$
from \eqref{Ann eq a1} and \eqref{Ann eq a3} then
\begin{equation}\label{Ann eq a4}
   D(f)=\frac{(z-\al_1)^{m_1}\ldots (z-\al_s)^{m_s}}{(z-\beta_1)^{n(n_1+k)}\ldots (z-\beta_t)^{n(n_t+k)}}g(z),
\end{equation}
where $g$ is a polynomial and
\begin{equation}\label{Ann eq a5}
\operatorname{deg}(g) \leq \max\{(n-1)N+nkt,(n-1)M-nks+n \operatorname{deg}(g_1)\}
\end{equation}

   Since $D(f)(z)-z^m$ has exactly one zero at $z_0$ with multiplicity $l,$ we have
   \begin{equation}\label{Ann eq a6}
    D(f(z))=z^m+\frac{B(z-z_0)^l}{(z-\beta_1)^{n(n_1+k)}\ldots (z-\beta_t)^{n(n_t+k)}},
   \end{equation}
   where $B$   is a non-zero constant and $l$ is a positive integer. On differentiating   \eqref{Ann eq a4} and \eqref{Ann eq a6} $m+1$ times, we get
   \begin{equation}\label{Ann eq a7}
    ( D(f))^{(m+1)}=\frac{(z-\al_1)^{m_1-m-1}\ldots (z-\al_s)^{m_s-m-1}h_1(z)}{(z-\beta_1)^{n(n_1+k)+m+1}\ldots (z-\beta_t)^{n(n_t+k)+m+1}},
   \end{equation}
where $h_1$ is a polynomial with $\operatorname{deg}h_1\leq (m+1)(s+t-1)+ \operatorname{deg}(g)$. And
   \begin{equation}\label{Ann eq a8}
     ( D(f))^{(m+1)}=\frac{(z-z_0)^{l-m-1}h_2(z)}{(z-\beta_1)^{n(n_1+k)+m+1}\ldots (z-\beta_t)^{n(n_t+k)+m+1}},
   \end{equation}
where $h_2$ is a polynomial with  $\operatorname{deg}h_2\leq(m+1)t.$\\

Since $\al_i\neq z_0$ for $i=1, 2, \ldots, s,$ it follows from \eqref{Ann eq a7} and \eqref{Ann eq a8} that 
\begin{equation*}M-(m+1)s\leq \operatorname{deg}(h_2)\leq (m+1)t,\end{equation*}which implies that \begin{equation*} M\leq (m+1)(s+t) <(k+1)(s+t)\leq(k+1)\left(\frac{M}{2k+2}+\frac{N}{2k+3}\right).\end{equation*}
Hence we deduce that
  \begin{equation}\label{Ann eq a9}
M<N
\end{equation}

Now we discuss the following two subcases.\\

\underline{Case1.1.1} If $l\neq m+nN+nkt.$ It follows from \eqref{Ann eq a4} that
\begin{equation}\label{Ann eq p1}
nN+nkt\leq M+ \operatorname{deg}(g).
\end{equation}
If $\operatorname{deg}(g)\leq (n-1)N+nkt,$ we thus from \eqref{Ann eq a4} obtain that $nq+nkt\leq M+(n-1)N+nkt,$ which implies that $N\leq M<N$ by \eqref{Ann eq a9}. This is impossible.\\

If $\operatorname{deg}(g)\leq (n-1)M-nks+n \operatorname{deg}(g_1),$ since $\operatorname{deg}(g_1)\leq k(s+t-1),$ hence $nN+nkt\leq(n-1)M-nks+nk(s+t-1),$ we have $N\leq M-1<N-1$ by \eqref{Ann eq a9}. We thus arrive at a contradiction .\\

\underline{Case1.1.2} When $l=m+nN+nkt.$ It is obtained from \eqref{Ann eq a7} and \eqref{Ann eq a8} that 
\begin{equation}\label{Ann eq p2}
l-m-1\leq\operatorname{deg}(h_1)\leq (m+1)(s+t-1)+\operatorname{deg}(g).
\end{equation}  
If $\operatorname{deg}(g)\leq (n-1)N+nkt,$ we thus from \eqref{Ann eq p2} obtain that \begin{equation*}
 l\leq (m+1)(s+t)+\operatorname{deg}(g),
\end{equation*} which implies that      $$m+nN+nkt\leq (m+1)(s+t)+(n-1)N+nkt.$$  We have, $$N\leq(m+1)(s+t)<(k+1)\left(\frac{M}{2k+2}+\frac{N}{2k+3)}\right)<M,$$ which is  a contradiction.\\

If $$\operatorname{deg}(g)\leq (n-1)M-nks+n\operatorname{deg}(g_1) \leq (n-1)M-nks+nk(s+t-1).$$ By eq \eqref{Ann eq p2}  $$m+nN+nkt\leq(m+1)(s+t)+(n-1)M-nks+nk(s+t-1).$$ This gives that $$nN\leq(m+1)(s+t)+(n-1)M-nk<(k+1)\left(\frac{M}{2k+2}+\frac{N}{2k+3}\right)+(n-1)M-nk.$$ Which gives $N<M-1<N-1,$ this is a contradiction.\\

\underline{Case2.} Let $D(f)(z)-z^m$ has no zero . Then $f$ can not be a polynomial. Hence $f$ is non polynomial rational function. Now putting $l=0$ in \eqref{Ann eq a6} and proceeding as Case $1.1$ of lemma, we have a contradiction. 

\end{proof}

\section{Proof of Main Results} First we give the proof of Theorem \ref{thm 1}. 
\begin{proof}[Proof of Theorem \ref{thm 1}]  Since normality is a local property, we assume that
${D}=\Delta$. Suppose that $\fr$ is not normal
in $\Delta$. Then there exists at least one point $z_0$ such that
$\fr$ is not normal  at the point $z_0$ in $\Delta$. Without loss of
generality we may assume that $z_0=0$. 
 By Lemma \ref{lemma 1}, there
exist
\begin{enumerate}
\item{ a sequence of complex numbers $z_j \rightarrow 0$, $|z_j|<r<1$},
\item{ a sequence of functions $f_j\in \mathcal F$, }
\item{ a sequence of positive numbers $\rho_j \rightarrow 0$},
\end{enumerate}
such that $g_j(\zeta)=\rho_j^{-k}f_j(z_j+\rho_j\zeta) $ converges to a non-constant meromorphic function $g(\zeta)$ on $\C$ with $g^{\#}(\zeta)\leq g^{\#}(0)=1$. Moreover, $g$ is of order at most two. \\
We see that
\begin{equation}\label{eq3.2}
f_j(z_j+\rho_j\zeta)+a(z_j+\rho_j\zeta)\left(f_j^{(k)}(z_j+\rho_j\zeta) \right)^n-\al(z_j+\rho_j\zeta)\rightarrow a(0)\left(g^{(k)}(\zeta) \right)^n-\al(0).
\end{equation}
locally uniformly with respect to spherical metric on every compact subsets of $\C$ which contains no poles of $g$.\\

Clearly  $a(0)\left(g^{(k)}(\zeta) \right)^n-\al(0)\not\equiv 0$. Therefore by Lemma \ref{lemma 2} and Lemma \ref{lemma 3}, we know that $a(0)\left(g^{(k)}(\zeta) \right)^n-\al(0)$  has at least two distinct zeros. Now we claim that $a(0)\left(g^{(k)}(\zeta) \right)^n-\al(0)$ has only one zero.\\

 Contrary to this, let $a(0)\left(g^{(k)}(\zeta) \right)^n-\al(0)$ has two distinct zeros at $\zeta_0$ and $\zeta_1$. Now choose a small positive number $\delta$ such that $\Delta(\zeta_0, \delta)\cap \Delta(\zeta_1, \delta)=\emptyset $ and $a\left(g^{(k)}(\zeta) \right)^n-b$ has no other zeros in $\Delta(\zeta_0, \delta)\cup \Delta(\zeta_1, \delta)$. By Hurwitz's theorem, there exist two sequences  $\{\zeta_{j}\}\subset \Delta(\zeta_0,\delta), \{\zeta_{1_j}\}\subset\Delta(\zeta_{1},\delta)$ converging to $\zeta_{0}, \ \text{and}\  \zeta_{1}$ respectively and from \eqref{eq3.2}, for sufficiently large $j$, we have
\begin{align*}
  f_j(z_j+\rho_j\zeta_j)+a(z_j+\rho_j\zeta_j)\left(f_j^{(k)}(z_j+\rho_j\zeta_j) \right)^n-\al(z_j+\rho_j\zeta_j) &=0,\\
  f_j(z_j+\rho_j\zeta_{1_j})+a(z_j+\rho_j\zeta_{1_j})\left(f_j^{(k)}(z_j+\rho_j\zeta_{1_j}) \right)^n-\al(z_j+\rho_j\zeta_{1_j})&= 0.
  \end{align*}
Since $f+a\left(f^{(k)}\right)^n$ and $g+a\left(g^{(k)}\right)^n$ share $\al$ in $\Delta$, therefore for any positive integer $m$ we have
\begin{align*}
f_m(z_j+\rho_j\zeta_j)+a(z_j+\rho_j\zeta_j)\left(f_m^{(k)}(z_j+\rho_j\zeta_j) \right)^n-\al(z_j+\rho_j\zeta_j) &=0,\\
  f_m(z_j+\rho_j\zeta_{1_j})+a(z_j+\rho_j\zeta_{1_j})\left(f_m^{(k)}(z_j+\rho_j\zeta_{1_j}) \right)^n-\al(z_j+\rho_j\zeta_{1_j})&= 0.
\end{align*}
Fix $m$ and take $j\rightarrow \ity$, then we get $z_j+\rho_j\zeta_j\rightarrow 0, \ z_j+\rho_j\zeta_{1_j}\rightarrow 0$ and
\begin{equation}\label{eq3.3}
f_m(0)+ a(0)\left(f_m^{(k)}(0) \right)^n-\al(0)=0.
\end{equation}

Since the zeros  are isolated, so for large values of $j$  we have $z_j+\rho_j\zeta_{j}=0= z_j+\rho_j\zeta_{1_j}$. Hence
\begin{equation}
\zeta_j=-\frac{z_j}{\rho_j}, \qquad \zeta_{1_j}=-\frac{z_j}{\rho_j}.
\end{equation}
Which contradicts the fact that $\Delta(\zeta_0, \delta)\cap \Delta(\zeta_1, \delta)=\emptyset$.\end{proof}

\begin{proof}[Proof of Theorem \ref{thm 1'}] As in the proof of Theorem \ref{thm 1} assume $D=\Delta$ and $z_0=0$. \\

\underline{Case 1.} When $\alpha(0)\neq 0$, then there exists $r>0$ such that $\fr$ is normal in $|z|<r$ by Theorem \ref{thm 1}.

\underline{Case 2.} When $\alpha(0)=0$. So, we can write $\alpha(z)=z^m\beta(z),$ where $m$ is a positive integer and $\beta(z)$ is a holomorphic function in $D$ such that $\beta(0)\neq0$.  Again assuming  that $\fr$ is not normal at $0$, then by Zalcman's Lemma there exist
\begin{enumerate}
\item{ a sequence of complex numbers $z_j \rightarrow 0$, $|z_j|<r<1$},
\item{ a sequence of functions $f_j\in \mathcal F$, }
\item{ a sequence of positive numbers $\rho_j \rightarrow 0$},
\end{enumerate}
such that $g_j(\zeta)=\rho_j^{-\frac{nk}{n-1}}f_j(z_j+\rho_j\zeta)$ converges compactly to a non-constant meromorphic function $g(\zeta)$. Again, we  have two cases to consider:\\

\underline{Subcase 2.1.} If $z_j/\rho_j\rightarrow\infty$.  Then we consider the following family
\begin{equation}\label{Eq VIET1}
  \mathcal{G}:= \left\{G_j(\zeta)=z_j^{-\frac{nk}{n-1}}f_j(z_j(1+\zeta)): f_j\in \fr\right\}
\end{equation}
defined on $D$. From \eqref{Eq VIET1} we  obtain
\begin{equation*}
  G_j^{(k)}=z_j^{-\frac{k}{n-1}}f_j^{(k)}(z_j(1+\zeta)).
\end{equation*}
Let us define $D(f)(z)=f(z)+a\left(f^{(k)}\right)^n(z)$, then we have
\begin{align*}
  D(f_j)(z_j(1+\zeta)) & = f_j(z_j(1+\zeta))+a\left\{f_j^{(k)}(z_j(1+\zeta)\right\}^n \\
  & = z_j^{\frac{nk}{n-1}}G_j(\zeta)+z_j^{\frac{nk}{n-1}}a\left\{G_j^{(k)}(\zeta)\right\}^n \\
   & = z_j^{\frac{nk}{n-1}}D(G_j)(\zeta).
\end{align*}
Now, by the hypothesis for each pair $f^1, f^2$ in $\fr$,
\begin{equation*}
  (D(f^1)-\alpha)(z_j(1+\zeta))=0 \ \text{if and only if \ } (D(f^2)-\alpha)(z_j(1+\zeta))=0.
\end{equation*}
This gives that
\begin{equation*}
 z_j^{\frac{nk}{n-1}}D(G^1)(\zeta)-\alpha(z_j(1+\zeta))=0 \ \text{if and only if \ } z_j^{\frac{nk}{n-1}}D(G^2)(\zeta)-\alpha(z_j(1+\zeta)=0
\end{equation*}
This means
\begin{align*}
  D(G^1)(\zeta)&=z_j^{m-\frac{nk}{n-1}}(1+\zeta)^m\beta(z_j(1+\zeta))\\
  & \text{if and only if }\\
  D(G^2)(\zeta)&=z_j^{m-\frac{nk}{n-1}}(1+\zeta)^m\beta(z_j(1+\zeta)),
\end{align*}
Since $z_j^{m-\frac{nk}{n-1}}(1+\zeta)^m\beta(z_j(1+\zeta))\neq 0$ at the origin therefore by the previous case $\mathcal G$ is normal in $D$, hence there exists a subsequence  $\{G_j\}$ (after renumbering) in $\mathcal{G}$ such that $G_j\rightarrow G$, compactly in $D$.\\

Now, if $G(0)\neq 0, $ then we have
\begin{align*}
  g_j(\zeta) & = \rho_j^{-\frac{nk}{n-1}} f_j(z_j+\rho_j\zeta)=\left(\frac{z_j}{\rho_j}\right)^{\frac{nk}{n-1}}G_j\left(\frac{\rho_j}{z_j}\zeta\right)
   \end{align*}
   which converges to $\infty$ compactly on $\C$, which is a contradiction. Thus we must have $G(0)=0 \text{ and \ } G^{(2k+1)}(0)\neq \infty.$ \\

    And  for each $\zeta\in \C$ we have
    \begin{equation}\notag
      g_j^{(2k+1)}(\zeta) = \left(\frac{\rho_j}{z_j}\right)^{\frac{(n-2)k+n-1}{n-1}}G_j^{(2k+1)}\left(\frac{\rho_j}{z_j}\zeta\right)\rightarrow 0.
    \end{equation}
    This implies $g^{(2k+1)}(\zeta)\rightarrow 0,$  since all zeros of $g(\zeta)$ have multiplicity at least $2k+2$, so $g(\zeta)$ is
a constant

    \underline{Subcase 2.2.} If $z_j/\rho_j\rightarrow w_0$, where $w_0$ is a finite complex number. Then we see that
    \begin{equation}\notag
      H_j(\zeta)=\rho_j^{-\frac{nk}{n-1}} f_j(\rho_j\zeta)=g_j\left(\zeta-\frac{z_j}{\rho_j}\right)\rightarrow g(\zeta-w_0):= H(\zeta)
    \end{equation} compactly on $\C$.
    
    Also from Lemma \ref{lemma trans 1p} and Lemma \ref{lemma 2P} we have $D(H)(\zeta)-\zeta^m$ has at least two zeros. Now we proceed as in the proof of Theorem \ref{thm 1}.


  \end{proof}

\begin{proof}[Proof of Theorem \ref{thm 2}] We again assume that
${D}=\Delta$. Suppose that $\fr$ is not normal
in $\Delta$. Then there exists at least one point $z_0$ such that
$\fr$ is not normal  at the point $z_0$ in $\Delta$. Without loss of
generality we may assume that $z_0=0$. Then by Lemma \ref{lemma 1}, there
exist
\begin{enumerate}
\item{ a sequence of complex numbers $z_j \rightarrow 0$, $|z_j|<r<1$},
\item{ a sequence of functions $f_j\in \mathcal F$},
\item{ a sequence of positive numbers $\rho_j \rightarrow 0$},
\end{enumerate}
such that $g_j(\zeta)=\rho_j^{-k}f_j(z_j+\rho_j\zeta) $ converges to a non-constant meromorphic function $g(\zeta)$ on $\C$ with $g^{\#}(\zeta)\leq g^{\#}(0)=1$. Moreover, $g$ is of order at most two. \\
We see that
\begin{equation}\label{eq3.2a}
f_j(z_j+\rho_j\zeta)+a(z_j+\rho_j\zeta)\left(f_j^{(k)}(z_j+\rho_j\zeta) \right)^n-b\rightarrow a(0)\left(g^{(k)}(\zeta) \right)^n-b.
\end{equation}
locally uniformly with respect to spherical metric on every compact subsets of $\C$ which contains no poles of $g$.\\

 Now we claim that $a(0)\left(g^{(k)}(\zeta) \right)^n-b$ has at most one zero IM. Suppose on contrary, let $a(0)\left(g^{(k)}(\zeta) \right)^n-b$ has two distinct zeros at $\zeta_0$ and $\zeta_1$. Now choose a small positive number $\delta$ such that $\Delta(\zeta_0, \delta)\cap \Delta(\zeta_1, \delta)=\emptyset $ and $a(0)\left(g^{(k)}(\zeta) \right)^n-b$ has no other zeros in $\Delta(\zeta_0, \delta)\cup \Delta(\zeta_1, \delta)$. By Hurwitz's theorem, there exist two sequences  $\{\zeta_{j}\}\subset \Delta(\zeta_0,\delta), \{\zeta_{1_j}\}\subset\Delta(\zeta_{1},\delta)$ converging to $\zeta_{0}, \ \text{and}\  \zeta_{1}$ respectively and from \eqref{eq3.2}, for sufficiently large $j$, we have
\begin{align*}
  f_j(z_j+\rho_j\zeta_j)+a(z_j+\rho_j\zeta_j)\left(f_j^{(k)}(z_j+\rho_j\zeta_j) \right)^n-b &=0,\\
  f_j(z_j+\rho_j\zeta_{1_j})+a(z_j+\rho_j\zeta_{1_j})\left(f_j^{(k)}(z_j+\rho_j\zeta_{1_j}) \right)^n-b&= 0.
  \end{align*}

For large values of $j$ $z_j+\rho_j\zeta_j\in \Delta(\zeta_0, \delta)$ and $z_j+\rho_j\zeta_{1_j}\in\Delta(\zeta_1, \delta)$, so $f_j+a\left(f_j^{(k)}\right)^n-b$ has two distinct zeros, which contradicts the fact that $f_j+a\left(f_j^{(k)}\right)^n-b$ has at most one zero. But Lemma \ref{lemma 2} and Lemma \ref{lemma 3} confirm  the non-existence of such non-constant meromorphic function $g$. This contradiction shows that $\fr$ is normal in $\Delta$ and this proves the theorem. \end{proof}

 The proof of Theorem \ref{lemma 6} is same as the proof of Theorem E   with  little changes in the last lines. For completeness we give the proof of Theorem \ref{lemma 6}.\\

\begin{proof}[Proof of Theorem \ref{lemma 6}]   Suppose $f+a\left(f^{(k)} \right)^n$ assumes each value $b\in \C$ only finitely many times. This means \begin{equation}\label{eq2.34}
N\left(r, \frac{1}{f+a\left(f^{(k)} \right)^n-b}\right)=o(\log r) = S(r, f)
\end{equation}

Let us define \begin{equation}\label{eq2.21}F:=f+a\left(f^{(k)}\right)^n-b,\end{equation}   \begin{equation}\label{eq2.22}\phi:=\frac{F'}{F}, \end{equation} \begin{equation}\label{eq2.23} \psi:=n\frac{f^{(k+1)}}{f^{(k)}}-\frac{F'}{F}.\end{equation}
Now, we claim that $\phi\psi\not\equiv 0.$ If $\phi\equiv0$, then $F'\equiv 0$. We can deduce that $F\equiv  c$, where $c$ is finite complex number. From \eqref{eq2.21} and Lemma \ref{lem} we get that  $f$ must be a polynomial, which is a contradiction.\\

Next, if $\psi\equiv 0$, from \eqref{eq2.23}, we get
\begin{equation}\label{eq2.24}
c\left(f^{(k)}\right)^n = f +a\left(f^{(k)}\right)^n -b,
\end{equation}
where $c\in\C$. From \eqref{eq2.24} we get
\begin{equation}\label{eq2.25}
(a-c)\left(f^{(k)}\right)^n + f =b.
\end{equation}
If $a-c=0,$ we get that $f\equiv b$, which is contradiction. Otherwise, we conclude from \eqref{eq2.25} and Lemma \ref{lem} that $f$ must be a polynomial, which is a contradiction. \\

From \eqref{eq2.21}, we have $T(r, F)=O(T(r, f))$, thus from \eqref{eq2.22}
and \eqref{eq2.23}, we have \begin{equation}\label{eq2.26} m(r, \phi) =S(r, f) \ \text{and} \ m(r, \psi)= S(r, f).\end{equation}

From \eqref{eq2.34}, \eqref{eq2.22}, \eqref{eq2.23} and Nevanlinna's First Fundamental Theorem, we get
\begin{align}
N\left(r, \frac{1}{\phi}\right)&\leq m(r,\phi) + N(r, \phi) - m\left(r, \frac{1}{\phi}\right) + O(1)\notag\\
&\leq N(r,\phi) + S(r, f) \leq \overline{N}(r, f) + S(r,f), \label{eq2.27}\\
N\left(r, \frac{1}{\psi}\right)&\leq m(r,\psi) + N(r, \psi) - m\left(r, \frac{1}{\psi}\right) + O(1)\notag\\
&\leq N(r,\psi) + S(r, f) \leq \overline{N}\left(r, \frac{1}{f^{(k)}}\right) + S(r,f). \label{eq2.28}
\end{align}
Again, by \eqref{eq2.22} and \eqref{eq2.23}, we get
\begin{equation}\label{eq2.29}
(f-b)\phi-f'= a\left(f^{(k)}\right)^n \psi.
\end{equation}
 So, we get from \eqref{eq2.34}, \eqref{eq2.26} and \eqref{eq2.27}
 \begin{align}
 T(r, &(f-b)\phi -f')= T\left(r, (f-b)\left(\phi-\frac{f'}{f-b}\right)\right)\notag\\
 &\leq T(r, f-b) + T\left(r, \phi - \frac{f'}{f-b}\right) +S(r,f)\notag\\
 &\leq m(r,f-b)+N(r,f-b)+ m\left(r,\phi-\frac{f'}{f-b}\right)+ N\left(r, \phi-\frac{f'}{f-b}\right)+S(r,f)\label{eq2.30}\\
 &\leq m(r,f) + N(r, f) + m(r,\phi) + m(r,\frac{f'}{f-b}) + N\left(r, \phi-\frac{f'}{f-b}\right)+S(r,f)\notag\\
 & \leq T(r, f) + \overline{N}(r, f) + S(r, f)\notag.
 \end{align}

 It follows from  \eqref{eq2.26} - \eqref{eq2.30} that
 \begin{align}
 nT\left(r,f^{(k)} \right)&\leq T(r, \psi) + T(r, \phi(f-b)-f')+ S(r, f)\notag\\
 & \leq m(r, \psi) + N(r, \psi) + T(r, f) + \overline{N}(r,f) + N\left(r, \frac{1}{F}\right)+ S(r,f)\notag\\
 & \leq \overline{N}\left(r, \frac{1}{f^{(k)}}\right) + N\left(r, \frac{1}{F}\right) + m\left(r, \frac{1}{f}\right)+ N\left(r, \frac{1}{f}\right)+ \overline{N}(r,f) + S(r,f)\label{eq2.31}\\
 & \leq \overline{N}\left(r, \frac{1}{f^{(k)}}\right)+ 2N\left(r, \frac{1}{F}\right)+ m\left(r,\frac{1}{f^{(k)}}\right)+N\left(r, \frac{1}{f}\right)+\overline{N}(r, f) + S(r,f)\notag\\
 &\leq T\left(r, \frac{1}{f^{(k)}}\right) +2N\left(r, \frac{1}{F}\right)+ N\left(r, \frac{1}{f}\right)+\overline{N}(r, f) +S(r,f)\notag\\
 & \leq T\left(r, f^{(k)}\right) + 2N\left(r, \frac{1}{F}\right)+ N\left(r, \frac{1}{f}\right)+\overline{N}(r, f) + S(r, f).\notag
 \end{align}
 Therefore, we have from \eqref{eq2.34},
 \begin{equation}\label{eq2.32}
 (n-1)T\left(r, f^{(k)}\right)\leq N\left(r, \frac{1}{f}\right)+\overline{N}(r, f)+ S(r,f).
 \end{equation}
Also, we have
\begin{equation}\label{eq2.33}
(n-1)T\left(r, f^{(k)}\right)\geq (n-1)N\left(r, f^{(k)} \right)\geq (n-1)N(r,f)+ (n-1)\overline{N}(r,f).
\end{equation}

Thus by \eqref{eq2.32} and \eqref{eq2.33} we have
\begin{align*}
 (n-1)N(r,f)+ (n-1)\overline{N}(r,f)&\leq N\left(r, \frac{1}{f}\right)+\overline{N}(r, f)+ S(r,f)\\
& \leq T\left(r, f\right)+\overline{N}(r, f)+ S(r,f).
\end{align*}
This gives
\begin{equation}\notag
(n-2)N(r,f)+ (n-2)\overline{N}(r,f)\leq S(r,f).
\end{equation}
So, we get \begin{equation}\label{eq2.35}N(r, f)=S(r,f).\end{equation}
Therefore from \eqref{eq2.32}, we get
\begin{align*}
(n-1)T\left(r, f^{(k)}\right)&\leq N\left(r, \frac{1}{f} \right)+\overline{N}(r, f) + S(r, f),\\
&\leq T(r,f) + S(r,f).
\end{align*}
And this gives
\begin{equation}\notag
(n-2)T\left(r, f^{(k)}\right)+T\left(r, f^{(k)}\right)-T(r,f)\leq S(r, f).
\end{equation}
Now, by \eqref{eq2.35}, we get
\begin{equation}
T\left(r, f^{(k)}\right)\leq S(r, f),
\end{equation}
which is a contradiction. So $f+a\left(f^{(k)}\right)^n$ assumes each value $b\in \C$ infinitely often. \end{proof}


\begin{thebibliography}{00}

\bibitem{Ahl} L. V. Ahlfors, {Complex Analysis,} Third edition, \emph{McGraw-Hill}, 1979.
\bibitem{CTZY} W. Chen et al, \emph{Normal families of meromorphic functions concerning shared values}, J. Comp. Anal., vol. 2013, Article ID 107281, 7 pages, 2013.
\bibitem{dethloff} G. Dethloff, T. V. Tan, N. V. Thin, \emph{Normal criteria for  families of meromorphic functions,} J. Math. Anal. Appl. 411 (2014), 675--683.
\bibitem{FZ} M. L. Fang, L. Zalcman, \emph{On the value distribution of $f+a(f')^n$,} Sci. China Ser. A: Math. 51 (2008), no. 7, 1196--1202.
\bibitem{Hay} W. K. Hayman, {Meromorphic Functions,} \emph{Clarendon Press, Oxford}, 1964.
\bibitem{hin} J. D. Hinchliffe, \emph{On a result of Chuang related to Hayman's alternative,} Comput. Methods Funct. Theory 2 (2002), 293--297.
\bibitem{LFZ} C. L. Lei et al, \emph{Some normality criteria of meromorphic functions,} Acta Math. Sci. 33 B (2013), no. 6, 1667--1674.
    \bibitem{YLI} Y. T. Li, \emph{Normal families of zero-free meromorphic functions}, Abst. and Appl. Anal., Vol. 2012, Article ID 908123, 12 pages, 2012.

\bibitem{MS}E. Mues and N. Steinmetz, \emph{Meromorphe Functionen, die mit ihrer Ableitung Werte teilen,} Manuscripta Math. 29 (1979), 195--206.
\bibitem{Schiff} J. Schiff, { Normal Families,} \emph{Springer-Verlag, Berlin,} 1993.
\bibitem{Sch}W. Schwick, \emph{Sharing Values and Normality,} Arch. Math. 59 (1992), 50--54.
\bibitem{YMW} Y. M. Wang, \emph{On normality of meromorphic functions with multiple zeros and sharing values}, Turk. J. Math. 36 (2012), 263--271.
\bibitem{XWL} Y. Xu, F. Wu and L. Liao, \emph{Picard values and normal families of meromorphic functions,} Proc. Royal Soc. Edinburgh, 139 A, (2009), 1091--1099.
\bibitem{ccyang} C. C. Yang, H. X. Yi, {Uniqueness theory of meromorphic functions,} \emph{Science Press/ Kluwer Academic Publishers,} 2003.
\bibitem{Yang}L. Yang, {Value Distribution Theory,} \emph{Springer-Verlag, Berlin,} 1993.
\bibitem{Zalc 1}L. Zalcman, \emph{A heuristic principle in complex function theory,} The Amer. Math. Monthly, {82} (1975), 813--817.
\bibitem{Zalc}L. Zalcman, \emph{Normal families: new perspectives,} Bull. Amer. Math. Soc., {35}, no. 3  (1998), 215--230.

\end{thebibliography}
\end{document}